\documentclass[runningheads]{llncs}
\usepackage{graphicx}

\begin{document}

\title{Pawlak, Belnap and the magical number seven}

\author{Salvatore Greco\inst{1} and Roman S{\l}owi\'{n}ski \inst{2}} 

\institute{
Deaprtment of Economics and Business, University of Catania,\\
Corso Italia, 55, 95129 Catania, Italy \and
Institute of Computing Science, Pozna\'{n} University of Technology,\\
 60-965 Pozna\'{n},
 and Institute for Systems Research,\\
Polish Academy of Sciences, 01-447 Warsaw, Poland\
}

\maketitle

\begin{abstract}
We are considering the algebraic structure of the Pawlak-Brouwer-Zadeh lattice to distinguish vagueness due to imprecision from ambiguity due to coarseness. We show that a general class of  many-valued logics useful for reasoning about data emerges from this context. All these logics can be obtained from a very general seven-valued logic which, interestingly enough, corresponds to a reasoning system developed by Jaina philosophers four centuries BC. In particular, we show how the  celebrated Belnap four-valued logic can be obtained from the very general seven-valued logic based on the  Pawlak-Brouwer-Zadeh lattice.
\end{abstract}

\section{Introduction}

The Brouwer-Zadeh lattice \cite{Cattaneo_Nistico} was introduced as an algebraic structure to handle vagueness through representation of each concept $X$ of universe $U$ by  pair $(A,B)$, $A,B \subseteq U, A\cap B= \emptyset$, where $A$ is the necessity kernel (the set of objects from $U$ belonging to $X$ without any doubt) and $B$ is the non-possibility kernel (the set of objects from $U$ that for sure do not belong to $X$). In rough set theory \cite{pawlak82,pawlak91}, the Brouwer-Zadeh lattice can be seen as an abstract model  \cite{Cattaneo1997,Cattaneo2018,cc2004} representing each concept $X$  through pair of elements $(A,B)$ where $A$ is the lower approximation (interior) and $B$ is the complement of the upper approximation (exterior) of $X$. In \cite{IPMUGMS}, an extension of the  Brouwer-Zadeh lattice, called Pawlak-Brouwer-Zadeh lattice, has been proposed, where a new operator, called Pawlak operator, assigns the pair $(C,D)$ to each concept $X$ represented by pair $(A,B)$,  such that $C$ and $D$ represent the lower approximations of $A$ and $B$, respectively. The rough set theory of Pawlak operator has been discussed in \cite{gms2018}. In this paper, we reconsider the Pawlak-Brouwer-Zadeh lattice from Pawlak's perspective of reasoning about data \cite{pawlak91} and Belnap's idea of automated computer  reasoning \cite{Belnap76}. In particular, we demonstrate that within the Pawlak-Brouwer-Zadeh lattice, several many-valued logics naturally arise from a very general seven-valued logic, which is useful for modeling automated data reasoning. It is interesting to observe that this seven-valued logic corresponds to a reasoning system of argumentation proposed by Jaina philosophers four centuries BC \cite{Burch64,Priest2008}. Furthermore, the seven-valued logic is interesting from a cognitive psychology perspective. According to the influential and highly cited article 'The Magical Number Seven, Plus or Minus Two: Some Limits on Our Capacity for Processing Information' by Miller \cite{Miller1956}, it is suggested that individuals can effectively handle approximately seven stimuli simultaneously. This limit applies to both one-dimensional absolute judgment and short-term memory.

To give an intuition of the seven-valued logic and the other logics deriving from it, let us consider the following example. Consider a database containing data about the symptoms and related diagnosed diseases for a number of patients. We can imagine that for each considered disease, there are three possible types of diagnosis: 
\begin{itemize}
	\item the patient has the disease,
  \item  the patient does not have the disease,
 \item it cannot be said whether the patient has or not the disease (because, for instance, other tests have to be done).
\end{itemize}

On the other hand, for each patient $x$ there may be a number of patients in the database with the same symptoms. In this situation, it seems reasonable to analyze the database by classifying patients according to the diagnoses received by them and all other patients with the same symptoms. Proceeding in this way, there will be seven possible states of truth of the following proposition: ``A patient with the same symptoms as $x$ has the disease''. Let us name this proposition by $DISx$ and enumerate the possible cases: 
\begin{itemize}
	\item all the patients with the same symptoms as $x$ have the disease, so that, according to the available data, proposition $DISx$ is true,  
	\item some patients with the same symptoms as $x$ have the disease and for the others one cannot say if there is or not the disease, so that, according to the available data, proposition $DISx$  is sometimes true,  
	\item for all the patients with the same symptoms as $x$ one cannot say if there is or not the disease, so that, according to the available data, proposition $DISx$  is unknown,  
	\item some of the patients with the same symptoms as $x$ have the disease and the others have not the disease, so that, according to the available data, proposition $DISx$  is contradictory,  
	\item among the patients with the same symptoms as $x$ there are some with the disease, some without the diseases and some for which one does not know if there is or not the disease, so that, according to the available data, proposition $DISx$  is fully contradictory,  
	\item some patients with the same symptoms as $x$ do not have the disease and for the others one cannot say if there is or not the disease, so that, according to the available data, proposition $DISx$  is sometimes false,  
	\item all the patients with the same symptoms as $x$ have not the disease, so that, according to the available data, proposition $DISx$  is false.
\end{itemize}
The above seven situations are, of course, very detailed, so for practical reasons it could be more convenient to aggregate some of them.
For example, one could distinguish between the following two situations, the first of which would suggest some treatment, while in the second situation it would be more appropriate to wait for the possible appearance of some other symptoms: 
\begin{itemize}
	\item the case in which among the patients with the same symptoms there is some patient with the disease and there is no patient without the disease, that is the case in which  proposition $DISx$ is true or sometimes true, and    
\item the other cases. 
\end{itemize}
To support decisions about triage, a reasonable approach could be the following: 
\begin{itemize}
	\item if proposition $DISx$ is true or sometimes true, patient $x$ should be hospitalized,
	\item if proposition $DISx$ is false or sometimes false, patient $x$ should not be hospitalized,
	\item in all other cases, an expert doctor should be called to examine patient $x$ to make the decision about hospitalization.
\end{itemize}	
To diagnose a complex disease, the database could be used as follows:
\begin{itemize}
	\item if proposition $DISx$ is true or sometimes true, patient $x$ is diagnosed as having the disease,
	\item if proposition $DISx$ is unknown, some further medical tests are required for patient $x$,
	\item if proposition $DISx$ is contradictory or fully contradictory, an expert doctor should be called to examine patient $x$ to make the diagnosis,
\item if proposition $DISx$ is false or sometimes false, patient $x$ is not diagnosed as having the disease.
\end{itemize}
One can see that the basic case related to the general database query defines a seven-valued logic with truth values ``true'', ``sometimes true'', ``unknown'', ``contradictory'', ``fully contradictory'', ``sometimes false'' and ``false''. 
The reasoning about the treatment decision is based on a two-valued logic - apply or not the treatment - derived from the basic seven-valued logic. The envisaged protocol for triage gives an example of a three valued-logic  grounded again on the basic seven-valued logic. Finally, the diagnostic procedure proposes a possible four-valued logic that can also be obtained from the seven-valued logic.

Let us explain now how the basic seven-valued logic is related to the rough set concept and to the Pawlak-Brouwer-Zadeh lattice. First, observe that the set of all  patients in the data set having the same symptoms as patient $x$ represents the equivalence class $[x]_R$ with respect to the indiscernibility relation $R$ defined in terms of symptoms on universe $U$ which is the data set of patients. In other words, we consider indiscernible two patients $w$ and $z$ in the data set if $w$ has the same symptoms as $z$. Observe that each disease is represented in the data set by pair $(A,B)$ where $A$ is the set of patients from $U$ having the disease, and $B$ is the set of patients from $U$ not having the disease. The operators of the Pawlak-Brouwer-Zadeh lattice permit to define other interesting sets of patients from $U$. For example, applying the Brouwer negation to $(A,B)$ we get $(A,B)^{\approx}=(B,U-B)$, and consequently we obtain set $U-B$, i.e., the complement in $U$ of the set of  patients without the disease, or, in other words, the set of patients for whom there is the disease or it is unknown if there is the disease. In the same perspective, applying the Kleene negation and the Brouwer negation to $(A,B)$ we obtain $(A,B)^{- \approx}=(A,U-A)$,  so that we get set $U-A$, i.e., the complement in $U$ of the set of patients with the disease, or, in other words, the set of  patients for whom there is not the disease or it is unknown if there is or there is not the disease. Other interesting sets that can be obtained using the operators of the Pawlak-Brouwer-Zadeh lattice are $A \cup B$, i.e., the set of all patients for whom there is the disease or there is not the disease, and $U-A-B$, i.e., the set of all patients for whom it is unknown if there is or there is not the disease. Applying rough set theory, we can compute the lower and the upper approximation of the above mentioned sets of patients from $U$. For example, the set of patients $x$ from $U$ for whom all the other patients with the same symptoms have the disease constitutes the lower approximation of $A$ denoted by $\underline{R}A$. The set of all patients $x$ from $U$ for whom there is at least one patient with the same symptoms that has the disease constitutes, instead, the upper approximation of $A$ denoted by $\overline{R}A$. Within the Pawlak-Brouwer-Zadeh lattice, the rough approximation can be obtained through the application of the Pawlak operator that assigns pair $(\underline{R}A,\underline{R}B)$ to  pair $(A,B)$, that is $(A,B)^L=(\underline{R}A,\underline{R}B)$. Remembering that $\overline{R}A=U-\underline{R}(U-A)$, using the Pawlak operator and the other operators of the Pawlak-Brouwer-Zadeh lattice we can obtain the upper approximation of  set $A$ starting from pair $(A,B)$ as follows: $(A,B)^{- \approx L \approx}=(U-\underline{R}(U-A),\underline{R}(U-A))=(\overline{R}(A),\underline{R}(U-A))$. Analogously, starting from pair $(A,B)$, one can obtain all the following rough approximations: $\underline{R}A,\underline{R}B,\underline{R}(U-A-B),\underline{R}(U-A),\underline{R}(U-B),\underline{R}(A \cup B),\overline{R}A,\overline{R}B,\overline{R}(U-A-B)$. Using these rough approximations we can define all the truth values of the seven-valued logic as follows:
\begin{itemize}
	\item set of patients from $U$ for whom proposition $DISx$ is true: $\underline{R}A$,  
	\item set of patients from $U$ for whom proposition $DISx$  is sometimes true: \newline $\underline{R}(U-B)\cap\overline{R}A\cap\overline{R}(U-A-B)$,  
	\item set of patients from $U$ for whom  proposition $DISx$  is unknown: \newline $\underline{R}(U-A-B)$,  
	\item set of patients from $U$ for whom  proposition $DISx$  is contradictory: \newline $\underline{R}(A\cup B)\cap\overline{R}A \cap \overline{R}B$,  
	\item set of patients from $U$ for whom  proposition $DISx$  is fully contradictory: \newline $\overline{R}A \cap \overline{R}B \cap \overline{R}(U-A-B)$,  
	\item set of patients from $U$ for whom proposition $DISx$  is sometimes false: \newline $\underline{R}(U-A)\cap\overline{R}B\cap\overline{R}(U-A-B)$,   
	\item set of patients from $U$ for whom proposition $DISx$  is false:  $\underline{R}B$.
\end{itemize}
Using the truth values of the seven-valued logic and the operators of the Pawlak-Brouwer-Zadeh lattice, one can reconstruct all the other logics  derived from the aggregation of some of the seven truth values. For example, the truth values of the above mentioned two-valued logic for decisions about the treatment can be formulated as follows:
\begin{itemize}
	\item set of patients from $U$ for whom proposition $DISx$ is true or sometimes true: $\underline{R}A \cup (\underline{R}(U-B)\cap\overline{R}A\cap\overline{R}(U-A-B))=\underline{R}(U-B)\cap\overline{R}A$, and    
\item the other cases: $\underline{R}(U-A-B) \cup (\underline{R}(A\cup B)\cap\overline{R}A \cap \overline{R}B) \cup (\overline{R}A \cap \overline{R}B \cap \overline{R}(U-A-B)) \cup (\underline{R}(U-A)\cap\overline{R}B\cap\overline{R}(U-A-B)) \cup \underline{R}B=U-\underline{R}(U-B)=\overline{R}B \cup \underline{R}(U-A-B)$. 
\end{itemize}

The framework we have discussed is also related to three-way decision which is an interesting research line that has been deeply investigated in recent years (see \cite{Yao2021} for a survey). In fact, the three basic truth values - true, false, unknown - define a three-valued logic, and the derived seven-valued logic is based on the non-empty subsets of $\{$true, false, unknown$\}$.

 This paper presents all the above considerations formally and in detail. It is organized as follows. In the next section, we  recall the Pawlak-Brouwer-Zadeh lattice and the basic elements of  rough set theory. In the following section, we present the seven-valued logic and the other logics that can be derived from it. The last section collects conclusions. 

\section{The Pawlak-Brouwer-Zadeh distributive De Morgan lattices and indiscernibility-based rough set theory}

This section recalls the Pawlak-Brouwer-Zadeh distributive De Morgan lattices \cite{IPMUGMS} and shows  it as an abstract model of the classical rough set model based on indiscernibility.

A system $\left\langle \Sigma, \wedge, \vee,  ^\prime, ^{\sim}, 0, 1  \right\rangle$ is a quasi-Brouwer-Zadeh distributive lattice \cite{{Cattaneo_Nistico}} if the following properties (1)-(4) hold:
\begin{enumerate}
	\item[(1)] $\Sigma$ is a distributive lattice with respect to the join and the meet operations $\vee$, $\wedge$ whose induced partial order relation is 
\begin{center}
$a \le b$ iff $a=a \wedge b$ (equivalently $b=a \vee b$)
\end{center}
Moreover, it is required that $\Sigma$ is bounded by the least element 0 and the greatest element 1:
\begin{center}
$\forall a \in \Sigma, \quad 0 \le a \le 1$
\end{center}
	\item[(2)] The unary operation $^\prime: \Sigma \rightarrow \Sigma$ is a Kleene (also Zadeh or fuzzy) complementation. In other words, for arbitrary $a,b \in \Sigma$,
	
\begin{enumerate}
	\item[(K1)] $a^{\prime \prime}=a$,
	\item[(K2)] $(a \vee b)^\prime=a^\prime \wedge b^\prime$,
	\item[(K3)] $a \wedge a^\prime \le b \vee b^\prime$.
\end{enumerate}

	\item[(3)] The unary operation $^\sim$ $: \Sigma \rightarrow \Sigma$ is a Brouwer (or intuitionistic) complementation. In other words, for arbitrary $a,b \in \Sigma$,
	
\begin{enumerate}
	\item[(B1)] $a \wedge a^{\sim \sim} = a$,
	\item[(B2)] $(a \vee b)^\sim$ = $a^\sim$ $\wedge$ $b^\sim$,
	\item[(B3)] $a \wedge a^\sim =0$.
\end{enumerate}

\item[(4)]The two complementations are linked by the interconnection rule which must hold for arbitrary $a \in \Sigma$:
\vspace{0,1 cm}

(in) $a ^\sim$  $\le a^\prime$.
\end{enumerate}

A structure $\left\langle \Sigma, \wedge, \vee,  ^\prime, ^{\sim}, 0, 1  \right\rangle$ is a Brouwer-Zadeh distributive lattice if it is a quasi-Brouwer-Zadeh distributive lattice satisfying the stronger interconnection rule:
\vspace{0,1 cm}

(s-in) $a^{\sim \sim}=a^{\sim \prime}$.
\vspace{0,1 cm}

A Brouwer-Zadeh distributive lattice satisfying the $\vee$ De Morgan property: \vspace{0,2 cm}

(B2a) $(a\wedge b)^\sim=a^\sim \vee b^\sim$
\vspace{0,2 cm}

\noindent is called a De Morgan Brouwer-Zadeh distributive lattice.

%
%
%
%

An approximation operator, called Pawlak operator \cite{IPMUGMS}, on a De Morgan Brouwer-Zadeh distributive lattice is an  unary operation $^A: \Sigma \rightarrow \Sigma$ for which the following properties hold: for $a,b \in \Sigma$
\begin{itemize}
	\item[A1)] $a^{A\prime}=a^{\prime A}$;
	\item[A2)] $a\le b$ implies $b^{A\sim}\le a^{A\sim}$;
	\item[A3)] $a^{A\sim} \le a^\sim$ ;
	\item[A4)] $0^A=0$;
	\item[A5)] $a^\sim=b^\sim$ implies $a^A\wedge b^A=(a \wedge b)^A$;
	\item[A6)] $a^A\vee b^A \le(a \vee b)^A$;
	\item[A7)] $a^{AA}=a^A$;
	\item[A8)] $a^{A \sim A}=a^{A \sim}$;
	\item[A9)] $(a^A \wedge b^A)^A=a^A \wedge b^A$.
\end{itemize}

\subsection{Pawlak-Brouwer-Zadeh lattices and rough set theory}

A \textit{knowledge base} $K=(U,R)$ is a relational system where $U\neq \emptyset$ is a finite set called the \textit{universe} and $R$ is an equivalence relation on $U$. For any $x \in U$, $[x]_R$ is its equivalence class. The quotient set $U/R$ is composed of all the equivalence classes of $R$ on $U$. Given the knowledge base $K=(U,R)$, one can associate the two subsets ${\underline R} X$ and ${\overline R} X$ to each subset $X \subseteq U$:
$${\underline R} X=\{x \in U: [x]_R \subseteq X\},$$
$${\overline R} X=\{x \in U: [x]_R \cap X \neq \emptyset\}.$$
${\underline R} X$ and ${\overline R} X$ are called the lower and the upper approximation of $X$, respectively.

Let us consider the set of all pairs $\left\langle A,B \right\rangle$ such that $A,B \subseteq U$ and $A \cap B = \emptyset$. We denote by $3^U$ the set of these pairs, i.e.,
$$3^U=\{\left\langle A,B \right\rangle: A,B \subseteq U \mbox{ and } A \cap B = \emptyset\}.$$  

Given a knowledge base $K=(U,R)$, we can define an unary operator \mbox{$^L: 3^U\rightarrow 3^U$}, as follows: for any $\left\langle A,B \right\rangle \in 3^U$ 
$$\left\langle A,B \right\rangle^L=\left\langle {\underline R}A, {\underline R} B \right\rangle.$$
Let us consider the following operations on $3^U$: 
\begin{center}
$\left\langle A,B\right\rangle \sqcap \left\langle C,D\right\rangle=\left\langle A \cap C, B \cup D\right\rangle$,
\end{center}
\begin{center}
$\left\langle A,B\right\rangle \sqcup \left\langle C,D\right\rangle=\left\langle A \cup C, B \cap D\right\rangle$,
\end{center}
\begin{center}
$\left\langle A,B\right\rangle^-=\left\langle B,A\right\rangle$, 
\end{center}
\begin{center}
$\left\langle A,B\right\rangle^{\approx}=\left\langle B,U-B\right\rangle$.
\end{center}

The following results hold \cite{IPMUGMS}. 
\\

\textbf{Proposition 1.} The structure $\left\langle 3^U, \sqcap, \sqcup, ^{-},^{\approx},^L,\left\langle \emptyset, U \right\rangle, \left\langle U, \emptyset \right\rangle \right\rangle$ is a Pawlak-Brouwer-Zadeh lattice. $\qed$
\\

\textbf{Proposition 2.} For every Pawlak-Brouwer-Zadeh lattice ${\cal L}_{PBZ}= \linebreak \left\langle \Sigma, \wedge, \vee,  ^\prime, ^{\sim}, ^A, 0, 1  \right\rangle$, satisfying the condition 
\begin{itemize}
	\item[(P)] there exists $c\in\Sigma$ for which $c=c^\prime$,
\end{itemize}
there is a knowledge  base $K=(U,R)$ such that the structure 
$$RS_{PBZ}(U,R)=\left\langle 3^U, \sqcap, \sqcup, ^{-},^{\approx},^L,\left\langle \emptyset, U \right\rangle, \left\langle U, \emptyset \right\rangle \right\rangle$$
is isomorphic to ${\cal L}_{PBZ}$.

\section{The seven-valued logic of the Pawlak-Brouwer-Zadeh lattice}

In the following, we shall identify a set $S \subseteq U$ with the pair $\left\langle S,U-S \right\rangle$. Given a knowledge  base $K=(U,R)$, for each pair $a=\left\langle A,B \right\rangle, A,B \subseteq U, A\cap B = \emptyset$, the following sets can be considered:
\begin{itemize}
	\item \textit{the true part} of $\left\langle A,B \right\rangle$:
	
${\mathbf T}(A,B)=\\
=\{x \in U: [x]_R \subseteq A\}=\underline{R}A	=\left\langle A,B \right\rangle^{L - \approx} = a^{A \prime \sim}$, 
	\item \textit{the sometimes true part} of $\left\langle A,B \right\rangle$: 
	
	${\mathbf sT}(A,B)=\\
	=\{x \in U: [x]_R \subseteq (U -B), [x]_R \cap A \neq \emptyset \mbox{ and } [x]_R \cap (U-A-B) \neq \emptyset \}=\\=\underline{R} (U-B) \cap \overline{R}A \cap \overline{R}(U-A-B)	=\\=\left\langle A,B \right\rangle^{ \approx  L  \approx} \sqcap \left(\left\langle A,B \right\rangle^{ \approx -} \sqcap \left\langle A,B \right\rangle^{ - \approx -} \right)^{L \approx -} \sqcap \left\langle A,B \right\rangle^{- \approx L \approx -} =\\
	= a^{\sim A  \sim} \wedge \left(a^{\sim \prime} \wedge a^{\prime \sim \prime}\right)^{A \sim \prime} \wedge a^{\prime \sim A \sim \prime}$, 
	\item \textit{the unknown part} of $\left\langle A,B \right\rangle$: 
	
	${\mathbf U}(A,B)=\\
	=\{x \in U: [x]_R \subseteq U-A-B \}=\underline{R} (U-A-B)	= \left(\left\langle A,B \right\rangle^{ \approx -} \sqcap \left\langle A,B \right\rangle^{ - \approx -} \right)^{L  - \approx} =\\
	=  \left(a^{\sim \prime} \wedge a^{\prime \sim \prime}\right)^{A  \prime \sim}$,
	\item \textit{the contradictory part} of $\left\langle A,B \right\rangle$: 
	
	${\mathbf K}(A,B)=\\
	=\{x \in U: [x]_R \subseteq A\cup B, [x]_R \cap B \neq \emptyset \mbox{ and } [x]_R \cap A \neq \emptyset   \}=\\ =\underline{R} (A \cup B) \cap \overline{R}A	\cap \overline{R}B= \\=\left(\left\langle A,B \right\rangle \sqcup \left\langle A,B \right\rangle^- \right)^{L  - \approx}  \sqcap \left\langle A,B \right\rangle^{- \approx L \approx -} \sqcap \left\langle A,B \right\rangle^{ \approx L \approx -}=\\
	=  \left(a \vee a^\prime \right)^{A  \prime \sim} \wedge a^{\prime \sim A \sim \prime} \wedge a^{ \sim A \sim \prime}$, 

\item \textit{the fully contradictory part} of $\left\langle A,B \right\rangle$: 
	
	${\mathbf fK}(A,B)=\\
	=\{x \in U:  [x]_R \cap A \neq \emptyset,  [x]_R \cap B \neq \emptyset \mbox{ and } [x]_R\cap (U-A -B) \neq \emptyset \}=\\ = \overline{R}A	\cap \overline{R}B \cap \overline{R}(U-A-B)=\\ 
= \left\langle A,B \right\rangle^{- \approx L \approx -} \sqcap \left\langle A,B \right\rangle^{ \approx L \approx -} \sqcap \left(\left\langle A,B \right\rangle^{\approx -} \wedge \left\langle A,B \right\rangle^{- \approx -}\right)^{L \approx -}  =\\
	=   a^{\prime \sim A \sim \prime} \wedge a^{ \sim A \sim \prime} \wedge \left(a^{\sim \prime}  \wedge a^{\prime \sim \prime} \right)^{A  \sim \prime}$,

\item \textit{the sometimes false part} of $\left\langle A,B \right\rangle$: 
	
	${\mathbf sF}(A,B)=\\
	=\{x \in U: [x]_R \subseteq (U-A), [x]_R \cap B \neq \emptyset \mbox{ and } [x]_R \cap (U-A-B) \neq \emptyset \}=\underline{R} (U-A) \cap \overline{R}B \cap \overline{R}(U-A-B)	=\\=\left\langle A,B \right\rangle^{ - \approx  L  \approx} \sqcap \left(\left\langle A,B \right\rangle^{ \approx -} \sqcap \left\langle A,B \right\rangle^{ - \approx -} \right)^{L \approx -} \sqcap \left\langle A,B \right\rangle^{- \approx L \approx -} =\\
	= a^{\prime \sim A  \sim} \wedge \left(a^{\sim \prime} \wedge a^{\prime \sim \prime}\right)^{A \sim \prime} \wedge a^{\prime \sim A \sim \prime}$,
\item \textit{the false part} of $\left\langle A,B \right\rangle$:
	
${\mathbf F}(A,B)=\\
=\{x \in U: [x]_R \subseteq B\}=\underline{R}B	=\left\langle A,B \right\rangle^{L \approx} = a^{A \sim}$. 
\end{itemize}
The truth values of the seven-valued logic can be represented by the lattice in Figure 1.
\begin{figure}[h]
	\centering
	\includegraphics[width=4cm]{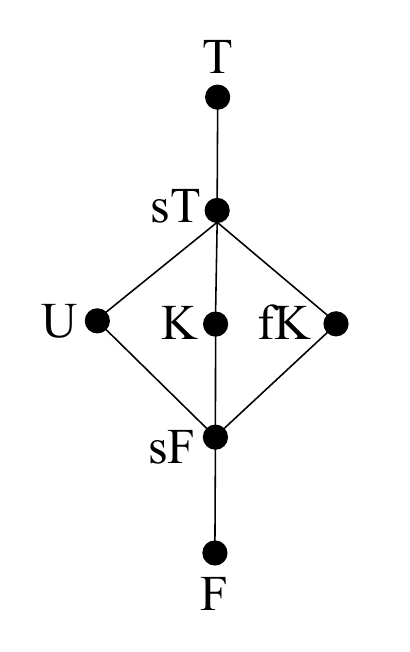}
	\caption{ \ Seven-valued logic truth value lattice }
	\label{figure 1}
\end{figure}

Let us remark that truth value operators of the seven-valued logic can be characterized in terms of upper approximations $\overline{R}A, \overline{R}B$ and $\overline{R}(U-A-B)$ as follows: for all $x \in U$
\begin{itemize}
	\item $x \in {\mathbf T}(A,B)=x \in \overline{R}A \wedge x \notin \overline{R}B  \wedge x \notin \overline{R}(U-A-B)$,
\item $x \in {\mathbf sT}(A,B)=x \in \overline{R}A \wedge x \notin \overline{R}B  \wedge  x \in \overline{R}(U-A-B)$,
\item $x \in {\mathbf U}(A,B)=x \notin \overline{R}A \wedge x \notin \overline{R}B  \wedge x \in \overline{R}(U-A-B)$,
\item $x \in {\mathbf K}(A,B)=x \in \overline{R}A \wedge x \in \overline{R}B  \wedge x \notin \overline{R}(U-A-B)$,
\item $x \in {\mathbf fK}(A,B)=x \in \overline{R}A \wedge x \in \overline{R}B  \wedge  x \in \overline{R}(U-A-B)$,
\item $x \in {\mathbf sF}(A,B)=x \notin \overline{R}A \wedge x \in \overline{R}B  \wedge x \in \overline{R}(U-A-B)$,
\item $x \in {\mathbf F}(A,B)=x \notin \overline{R}A \wedge x \in \overline{R}B  \wedge x \notin \overline{R}(U-A-B)$.
\end{itemize}
  
Observe that for all $A,B \subseteq U, A\cap B = \emptyset$,  
\[
{\mathbf T}(A,B) \cup {\mathbf sT}(A,B) \cup {\mathbf U}(A,B) \cup {\mathbf K}(A,B) \cup {\mathbf fK}(A,B) \cup {\mathbf sF}(A,B) \cup {\mathbf F}(A,B) = U
\]
and for all pairs $({\mathbf O_1},{\mathbf O_2})$ with ${\mathbf O_1},{\mathbf O_2} \in \mathcal{O}$$=\{{\mathbf T}, {\mathbf sT},  {\mathbf U} , {\mathbf K}, {\mathbf fK},  {\mathbf sF},  {\mathbf F}\}$ with ${\mathbf O_1} \neq {\mathbf O_2}$
\[
{\mathbf O_1}(A,B) \cap {\mathbf O_2}(A,B) = \emptyset.
\] 

Several aggregations of the truth value operators in $\mathcal{O}$ are interesting. Among them, the following upward and downward aggregations are particularly interesting:
\begin{itemize}
	\item \textit{the at least true part} of $\left\langle A,B \right\rangle$: 
	${\mathbf T^\uparrow}(A,B)={\mathbf T}(A,B)$,
\item \textit{the at least sometimes true part} of $\left\langle A,B \right\rangle$: 
	
	${\mathbf sT^\uparrow}(A,B)={\mathbf sT}(A,B) \cup {\mathbf T}(A,B) =\\=\{x \in U: [x]_R \subseteq U-B \mbox{ and } [x]_R \cap A \neq \emptyset\}=\underline{R} (U-B) \cap \overline{R}A=\\=	\left\langle A,B \right\rangle^{ \approx  L  \approx} \sqcap \left\langle A,B \right\rangle^{ - \approx  L  \approx -} = a^{\sim A  \sim} \wedge a^{\prime \sim A \sim \prime}$,
	
\item \textit{the at least unknown part} of $\left\langle A,B \right\rangle$: 
	
	${\mathbf U^\uparrow}(A,B)={\mathbf U}(A,B) \cup {\mathbf sT}(A,B) \cup {\mathbf T}(A,B)\\
	=\{x \in U: [x]_R \subseteq U-B \}=\underline{R} (U-B)	=\\
=\left\langle A,B \right\rangle^{ \approx L \approx }  	=   a^{ \sim A \sim}$,
	
\item \textit{the at least contradictory part} of $\left\langle A,B \right\rangle$: 
	
	${\mathbf K^\uparrow}(A,B)={\mathbf K}(A,B) \cup {\mathbf sT}(A,B) \cup {\mathbf T}(A,B)=\\
	=\{x \in U: \left([x]_R \subseteq (A \cup B) \mbox{ or } [x]_R \subseteq (U - B)\right) \mbox{ and } [x]_R \cap A \neq \emptyset\} =\\=\left(\underline{R} (A \cup B) \cup \underline{R} (U-B)\right) \cap \overline{R}A 	=\\= \left((\left\langle A,B \right\rangle \sqcup \left\langle A,B \right\rangle^-)^{ L - \approx} \sqcup \left\langle A,B \right\rangle^{ \approx L \approx}\right) \sqcap \left\langle A,B \right\rangle^{- \approx L \approx -}=\\ 	=  \left((a \vee a^\prime)^{ A \prime \sim} \vee a^{\sim A \sim}\right) \vee a^{\prime \sim A \sim \prime}$,

\item \textit{the at least fully contradictory part} of $\left\langle A,B \right\rangle$: 
	
	${\mathbf fK^\uparrow}(A,B)={\mathbf fK}(A,B)\cup{\mathbf sT}(A,B)\cup {\mathbf T}(A,B)\\
	=\{x \in U:  [x]_R \subseteq A \mbox{ or } ([x]_R \cap A \neq \emptyset \mbox{ and } [x]_R\cap (U-A -B) \neq \emptyset \}=\\ = \underline{R}A	\cup (\overline{R}A \cap \overline{R}(U-A-B)=\\=  \left\langle A,B \right\rangle^{ L - \approx} \sqcup (\left\langle A,B \right\rangle^{ - \approx L \approx -} \sqcap \left(\left\langle A,B \right\rangle^{\approx -} \sqcap \left\langle A,B \right\rangle^{- \approx -}\right)^{L \approx -}  =\\
	=   a^{A \prime \sim} \vee \left(a^{\prime  \sim A \sim \prime} \wedge (a^{\sim \prime}  \wedge a^{\prime \sim \prime})^{A \sim \prime}\right)$,

\item \textit{the at least sometimes false part} of $\left\langle A,B \right\rangle$: 
	
	${\mathbf sF^\uparrow}(A,B)={\mathbf sF}(A,B)\cup {\mathbf U}(A,B) \cup {\mathbf K}(A,B)\cup {\mathbf fK}(A,B)\cup{\mathbf sT}(A,B)\cup {\mathbf T}(A,B)\\
	=\{x \in U: [x]_R \cap (U-B) \neq \emptyset \}=\overline{R} (U-B) =	\left\langle A,B \right\rangle^{   L  \approx -} =  a^{A \sim \prime}$,
\item \textit{the at least false part} of $\left\langle A,B \right\rangle$:
	
${\mathbf F^\uparrow}(A,B)= U = \left\langle U, \emptyset \right\rangle=1$,

	\item \textit{the at most false part} of $\left\langle A,B \right\rangle$: 
	${\mathbf F^\downarrow}(A,B)={\mathbf F}(A,B)$,
\item \textit{the at most sometimes false part} of $\left\langle A,B \right\rangle$: 
	
	${\mathbf sT^\downarrow}(A,B)={\mathbf sF}(A,B) \cup {\mathbf F}(A,B) =\\=\{x \in U: [x]_R \subseteq U-A \mbox{ and } [x]_R \cap B \neq \emptyset\}=\underline{R} (U-A) \cap \overline{R}F=\\=	\left\langle A,B \right\rangle^{ - \approx  L  \approx} \sqcap \left\langle A,B \right\rangle^{  \approx  L  \approx -} = a^{\prime \sim A  \sim} \wedge a^{ \sim A \sim \prime}$,
	
\item \textit{the at most unknown part} of $\left\langle A,B \right\rangle$: 
	
	${\mathbf U^\downarrow}(A,B)={\mathbf U}(A,B) \cup {\mathbf sF}(A,B) \cup {\mathbf F}(A,B)\\
	=\{x \in U: [x]_R \subseteq U-A \}=\underline{R} (U-A)	=\\
=\left\langle A,B \right\rangle^{ - \approx L \approx }  	=   a^{ \prime \sim A \sim}$,
	
\item \textit{the at most contradictory part} of $\left\langle A,B \right\rangle$: 
	
	${\mathbf K^\downarrow}(A,B)={\mathbf K}(A,B) \cup {\mathbf sF}(A,B) \cup {\mathbf F}(A,B)=\\
	=\{x \in U: \left([x]_R \subseteq (A \cup B) \mbox{ or } [x]_R \subseteq (U - A)\right) \mbox{ and } [x]_R \cap B \neq \emptyset\} =\\=\left(\underline{R} (A \cup B) \cup \underline{R} (U-A)\right) \cap \overline{R}B 	=\\= \left((\left\langle A,B \right\rangle \sqcup \left\langle A,B \right\rangle^-)^{ L - \approx} \sqcup \left\langle A,B \right\rangle^{- \approx L \approx}\right) \sqcap \left\langle A,B \right\rangle^{ \approx L \approx -}=\\ 	=  \left((a \vee a^\prime)^{ A \prime \sim} \vee a^{\prime \sim A \sim}\right) \vee a^{ \sim A \sim \prime}$,

\item \textit{the at most fully contradictory part} of $\left\langle A,B \right\rangle$: 
	
	${\mathbf fK^\downarrow}(A,B)={\mathbf fK}(A,B)\cup{\mathbf sF}(A,B)\cup {\mathbf F}(A,B)\\
	=\{x \in U:  [x]_R \subseteq B \mbox{ or } ([x]_R \cap B \neq \emptyset \mbox{ and } [x]_R\cap (U-A -B) \neq \emptyset \}=\\ = \underline{R}B	\cup (\overline{R}B \cap \overline{R}(U-A-B)=\\=  \left\langle A,B \right\rangle^{ L \approx} \sqcup (\left\langle A,B \right\rangle^{  \approx L \approx -} \sqcap \left(\left\langle A,B \right\rangle^{\approx -} \sqcap \left\langle A,B \right\rangle^{- \approx -}\right)^{L \approx -}  =\\
	=   a^{A \sim} \vee \left(a^{ \sim A \sim \prime} \wedge (a^{\sim \prime}  \wedge a^{\prime \sim \prime})^{A \sim \prime}\right)$,

\item \textit{the at most sometimes true part} of $\left\langle A,B \right\rangle$: 
	
	${\mathbf sT^\downarrow}(A,B)={\mathbf sT}(A,B)\cup {\mathbf U}(A,B) \cup {\mathbf K}(A,B)\cup {\mathbf fK}(A,B)\cup{\mathbf sF}(A,B)\cup {\mathbf F}(A,B)\\
	=\{x \in U: [x]_R \cap (U-A) \neq \emptyset \}=\overline{R} (U-A) =	\left\langle A,B \right\rangle^{   L  - \approx } =  a^{A \prime \sim}$,
\item \textit{the at most true part} of $\left\langle A,B \right\rangle$:
	
${\mathbf T^\downarrow}(A,B)= U = \left\langle U, \emptyset \right\rangle=1$.
 
\end{itemize}

Let us consider the set of the upward truth value operators 

$\mathcal{O}$$^\uparrow=\{\mathbf{T^\uparrow},\mathbf{sT^\uparrow},\mathbf{U^\uparrow},\mathbf{K^\uparrow},\mathbf{fK^\uparrow},\mathbf{sF^\uparrow},\mathbf{F^\uparrow}\}$

and the set of the downward truth value operators

$\mathcal{O}$$^\downarrow=\{\mathbf{T^\downarrow},\mathbf{sT^\downarrow},\mathbf{U^\downarrow},\mathbf{K^\downarrow},\mathbf{fK^\downarrow},\mathbf{sF^\downarrow},\mathbf{F^\downarrow}\}.$

On the basis of $\mathcal{O}$$^\uparrow$ and $\mathcal{O}$$^\downarrow$ one $n$-valued logic with respect to the knowledge base $K=(U,R)$ is defined by the set of truth value operators $(O_1, \ldots O_n)$ such that, for all $A,B \subseteq U, A \cap B = \emptyset$, we have
\begin{itemize}
	\item for all $i=1, \ldots,n$
		\begin{itemize}
		\item either $O_i(A,B)=M_{1}(A,B) \cup \ldots \cup M_{k}(A,B)$ with $M_{1},\ldots,M_{k}\in$ $\mathcal{O}$$^\uparrow$ or
		\item $O_i(A,B)=M_{1}(A,B) \cup \ldots \cup M_{k}(A,B)$ with $M_{1},\ldots,M_{k}\in$ $\mathcal{O}$$^\downarrow$ or
		\item $O_i(A,B)=\bigg(M_{1}(A,B) \cup \ldots \cup M_{h}(A,B)\bigg)\cap \bigg(M_{h+1}(A,B) \cup \ldots \cup M_{k}(A,B)\bigg)$ with $M_{1},\ldots,M_{h}\in$ $\mathcal{O}$$^\uparrow$ and $M_{h+1},\ldots,M_{k}\in$ $\mathcal{O}$$^\downarrow$,
		\end{itemize}
		and 
	\item for all $O_i,O_j\in(O_1, \ldots O_n)$, $O_i(A,B)\cap O_j(A,B) = \emptyset$,
	\item $\bigcup_{i=1}^n O_i(A,B)=U$.		
\end{itemize}
  
Among the many posible logics that can be defined in this way, consider the Belnap four-valued logic \cite{Belnap77} defined by the set of truth value operators 
$$({\mathbf T}_{Belnap},{\mathbf U}_{Belnap},{\mathbf K}_{Belnap},{\mathbf F}_{Belnap})$$
 with
\begin{itemize}
	\item ${\mathbf T}_{Belnap}(A,B)={\mathbf sT^\uparrow}(A,B)={\mathbf T}(A,B) \cup {\mathbf sT}(A,B)=\overline{R}A \cap \underline{R}(U-B)$,
	\item ${\mathbf U}_{Belnap}(A,B)={\mathbf U}(A,B)={\mathbf U^\uparrow}(A,B) \cap {\mathbf U^\downarrow}(A,B)=\underline{R}(U-A-B)$,
	\item ${\mathbf K}_{Belnap}(A,B)=({\mathbf K^\uparrow}(A,B) \cup {\mathbf fK^\uparrow}(A,B))\cap ({\mathbf K^\downarrow}(A,B) \cup {\mathbf fK^\downarrow}(A,B))={\mathbf K}(A,B) \cup {\mathbf fK}(A,B)=\overline{R}A \cap \overline{R}B$,
	\item ${\mathbf F}_{Belnap}(A,B)={\mathbf sF^\uparrow}(A,B)={\mathbf F}(A,B) \cup {\mathbf sF}(A,B)=\overline{R}B \cap \underline{R}(U-A)$.
\end{itemize}
The above Belnap four-valued logic can be interpreted in terms of rough approximations as follows. Consider the knowledge base $K=(U,R)$, $x \in U$, and the concept $\left\langle A,B \right\rangle \in 3^U$. For all $x \in U$ we have
\begin{itemize}
	\item  $[x]_R \cap A \neq \emptyset$, i.e., $x \in \overline{R}A$, is an argument for truth, 
	\item  $[x]_R \cap B \neq \emptyset$, i.e., $x \in \overline{R}B$, is an argument for falsehood.
	\end{itemize}
Consequently,
\begin{itemize}
	\item  if $x \in \overline{R}A$ and $x \notin \overline{R}B$ (which is equivalent to $x \in \underline{R}(U-B)$), there are arguments for truth and there  are no arguments for falsehood, so that $x \in {\mathbf T}_{Belnap}(A,B)$,
	\item  if $x \notin \overline{R}A$ and $x \notin \overline{R}B$ (which is equivalent to $x \in \underline{R}(U-A-B)$), there are no arguments for truth and there  are no arguments for falsehood, so that $x \in {\mathbf U}_{Belnap}(A,B)$,
\item  if $x \in \overline{R}A$ and $x \in \overline{R}B$ , there are arguments for truth and there are arguments for falsehood, so that $x \in {\mathbf K}_{Belnap}(A,B)$,
	\item  if $x \notin \overline{R}A$ and $x \in \overline{R}B$ (which is equivalent to $x \in \underline{R}(U-A)$), there are no arguments for truth and there  are arguments for falsehood, so that $x \in {\mathbf F}_{Belnap}(A,B)$.
	
	\end{itemize}
Comparing the truth value operators of seven-valued logic, expressed in terms of upper approximations $\overline{R}A, \overline{R}B$ and $\overline{R}(U-A-B)$, with the Belnap's truth value operators of the four-valued logic reveals that the latter is different because it does not consider \mbox{$\overline{R}(U-A-B)$}.

\section{Conclusions}

The seven-valued logic considered in this paper naturally arises within the rough set framework, allowing to distinguish vagueness due to imprecision from ambiguity due to coarseness. 
We discussed the usefulness of this seven-valued logic  for reasoning about data. We showed that the Pawlak-Brouwer-Zadeh lattice is the proper algebraic structure for this seven-valued logic. We proposed also a general framework permitting to obtain other interesting many-valued logics by aggregation of truth value operators in the basic seven-valued logic. We plan to continue our research in this direction, investigating typical rough set topics, such as calculation of reducts and rule induction. We intend also to extend the seven-valued logic to reasoning about ordered data using the dominance-based rough set approach \cite{{gms2001}} and the related bipolar Pawlak-Brouwer-Zadeh lattice \cite{gms2012,gms2017}. Another interesting line of research we want to pursue is related to investigation of connections with other algebra models for rough sets such as Nelson algebra, Heyting algebra, \L ukasiewicz algebra, Stone algebra and so on (see, e.g., chapter 12 in \cite{Polkowski}). We also propose to study the relations between the rough set approach we adopted to obtain the seven-valued logic and the other logics that can be derived from it, and the approach presented in \cite{Nakayama2018} that derives four-valued logic from variable precision rough set model \cite{Ziarko1993}. Finally, with respect to the seven-valued logic we proposed, we plan to investigate the tableau calculi for deduction systems, taking into account also soundness and completeness.

\subsubsection*{Acknowledgments.} Salvatore Greco wishes to acknowledge the support of the Ministero dell’Istruzione, dell’Universit\'{a} e della Ricerca (MIUR) - PRIN 2017, project ``Multiple Criteria Decision Analysis and Multiple Criteria Decision Theory'', grant 2017CY2NCA. The research of Roman  S\l owi\'nski
 was supported by the SBAD funding from the Polish Ministry of Education and Science. This research also contributes to the PNRR GRInS Project.

\end{document}